\documentclass[11pt]{article}
\newtheorem{theorem}{Theorem}
\usepackage{graphicx}
\usepackage{amssymb}

\title{A sign pattern that allows oppositely signed orthogonal matrices}
\author{Bryan L. Shader \and Chanyoung Lee Shader\thanks{Mathematics Department 3036, 1000 E. University Avenue, Laramie, WY 82071.  \{bshader,chan\}@uwyo.edu}}
\date{\today}                                           

\begin{document}
\maketitle
\begin{abstract}
We provide the first example of a sign pattern $S$ for which there exist orthogonal matrices $Q_1$ and  $Q_2$
with sign pattern $S$ such that $\det Q_1=1$ and $\det Q_2=-1$.  The existence of such matrices is raised by C.~Waters  in {``Sign Pattern Matrices That Allow Orthogonality''}, {\it Linear Algebra and Its Applications},  235:1--13 (1996). 
\end{abstract}



\section{The question}
The {\it sign} of a real number $a$ is denoted by $\mbox{sgn}(a)$ and given by
$$
\mbox{sgn}(a)=\left\{ \begin{array}{rl}-1 & \mbox{ if $a<0$}\\
0 & \mbox{ if $a=0$}\\
1 & \mbox{ if $a>0$.}\end{array} 
\right.
$$
The {\it sign pattern} of a real matrix $A$ is the matrix obtained from $A$ by replacing each of 
its entries by its sign.  An $n\times n$ matrix $S$ determines a class $\mathcal{Q}(S)$ 
consisting of all matrices with sign pattern $S$.

Several papers \cite{BBS, CS, F, JW,W} have studied the properties of sign patterns of orthogonal matrices.
The sign pattern $S$ {\it allows orthogonality} if there exists an orthogonal matrix $U$ in $\mathcal{Q}(S)$.
For example, 
$$
S=\left[ \begin{array}{rrr}
1 & -1 & 0 \\
1 & 1 & -1 \\
1 & 1 & 1 \end{array}
\right]
$$
allows orthogonality since 
$$
\left[\begin{array}{rrr}
\frac{1}{\sqrt{2}} & - \frac{1}{\sqrt{2}}&0 \\
\frac{1}{2} & \frac{1}{2} & -\frac{1}{\sqrt{2}}\\
\frac{1}{2} & \frac{1}{2} & \frac{1}{\sqrt{2}}
\end{array} \right]
$$
is an orthogonal matrix in $\mathcal{Q}(S)$.
Also, 
$$
T= 
\left[ \begin{array}{rrr}
1 & 1 & 0 \\
1 & 1 & -1 \\
1& 1 & 1 \end{array}
\right]
$$
does not allow orthogonality since the dot-product between first and second columns of 
each matrix in $\mathcal{Q}(T)$ is positive. 

An elementary fact is that the determinant of an orthogonal matrix is either $1$ or $-1$. 
In \cite{W}, C. Waters observed that for each sign pattern $S$ of order at most $4$,
either each orthogonal matrix in $\mathcal{Q}(S)$ has determinant $1$, or 
each orthogonal matrix in $\mathcal{Q}(S)$ has determinant $-1$.  He also noted 
that for the $n$ by $n$ sign pattern $S$ consisting of $-1$'s in positions $(k,k)$ $(k=2,\ldots, n)$ and $1$'s elsewhere
has the property that every orthogonal matrix with sign pattern $S$ has determinant $(-1)^{n-1}$. 

In this note
we provide the first example, to our knowledge, of a pair of orthogonal matrices with the same 
sign pattern whose determinants have opposite sign.  In a forthcoming paper, we will describe
the theory that led to the construction of these matrices and further results concerning sign patterns 
that allow orthogonality.

\section{The matrices}

Let
$$Q_1=\frac{1}{8}
\left[ \begin{array}{rrrrrrr}
5 & 3 & 2& -3& -2& \; 3&  2\\
-2 & 5& -3 & 2 & 3 &\; 3&  2\\
-3&  2 & 5 & 3& -3& \; 2& -2\\
 2& -3& -2&  5& -3& \; 2&  3\\
3 &-2 & 2&  2&  5& \; 3 &-3\\
-2& -2& -3& -3& -2&\;  5& -3\\
-3 &-3&  3& -2&  2& \; 2&  5\end{array}
\right]
$$
and
$$
Q_2=\frac{1}{20014} \left[ \begin{array}{rrrrrrr}
 9389  &  396 & 10197 &  -396 &-10197 &   396 & 10197\\
-10197 & 9389 & -396 &10197 &  396 &  396 & 10197\\
  -396&  10197 &  9389 &   396 &  -396&  10197& -10197\\
10197 &  -396 &-10197 &  9389 &  -396&  10197 &   396\\
396 & -10197 & 10197&  10197 &  9389&    396&   -396\\
-10197 &-10197 &  -396 &  -396& -10197 &  9389&   -396\\
-396 &  -396 &   396 &-10197 & 10197 & 10197 &  9389
\end{array}
\right]
.
$$
It is not difficult to verify (we used Sage \cite{S}) that both $Q_1$ and $Q_2$ are orthogonal,
$\det Q_1>0$, $\det Q_2<0$ and $Q_1$ and $Q_2$ have the same sign pattern.

\end{document}